\documentclass[11pt]{article}
\usepackage{amsfonts}
\usepackage{mathrsfs}
\usepackage{amssymb,amsmath} 
 \allowdisplaybreaks

\catcode`!=11
\let\!int\int \def\int{\displaystyle\!int}
\let\!lim\lim \def\lim{\displaystyle\!lim}
\let\!sum\sum \def\sum{\displaystyle\!sum}
\let\!sup\sup \def\sup{\displaystyle\!sup}
\let\!inf\inf \def\inf{\displaystyle\!inf}
\let\!cap\cap \def\cap{\displaystyle\!cap}
\let\!max\max \def\max{\displaystyle\!max}
\let\!min\min \def\min{\displaystyle\!min}
\let\!frac\frac \def\frac{\displaystyle\!frac}
\catcode`!=12

\let\oldsection\section
\renewcommand\section{\setcounter{equation}{0}\oldsection}

\allowdisplaybreaks
\def\pf{\it{Proof.}\rm\quad}

\def\R{\mathbb{R}}

\newtheorem{thm}{Theorem}[section]

\newtheorem{pro}{Proposition}[section]
\newtheorem{lem}{Lemma}[section]

\begin{document}
\title{\Large\bf  Asymptotic behavior of solutions toward the strong contact discontinuity for  compressible Navier-Stokes equations with Cauchy problem}
\author{{Tingting Zheng\ \footnote{Corresponding author email:
asting16@sohu.com(T.Zheng),}\footnote{
 This work was partially supported by
the Youth Natural Science Foundation of Fujian Province, China
(Grant No. 2017J05001) .}}\\[2mm]
{\small Computer and Message Science College, Fujian Agriculture and
Forest University,}\\{\small Fuzhou 350001, P. R. China}\\[2mm]
}

\date{}

\maketitle

\noindent{\bf Abstract.}In this paper, we consider the nonisentropic
ideal polytropic Navier-Stokes equations to the Cauchy problem.  The
asymptotic stability of contact discontinuity is established under
the condition that the initial perturbations are partly small but
the strength of contact discontinuity can be suitably large. With
this conditions,  the bounds of density and temperature can be
obtained from the complicated structure of Navier-Stokes equations.
The proofs are given by the elementary energy method.
 \\[2mm]
\noindent{\bf AMS Subject Classifications (2000).} 35B40, 35B45,
76N10,76N17
\\[2mm]
\noindent{\bf Keywords:}Cauchy problem, Compressible Navier-Stokes equations,
Strong contact discontinuity, Asymptotic stability .
\section{Introduction}
This paper is concerned with the one-dimensional compressible
viscous heat-conducting flows in the whole space
$\R=(-\infty,+\infty)$ , which is governed by the following initial
value problem in Eulerian coordinate $(\tilde{x},t)$:
\begin{equation}
\left\{
\begin{array}{lll}
\tilde\rho_t+(\tilde\rho \tilde u)_{\tilde x}=0,\quad (\tilde x,t)\in\R\times\R_+,\\[2mm]
(\tilde\rho \tilde u)_t+(\tilde\rho \tilde u^2+\tilde p)_{\tilde{x}}=\mu \tilde u_{\tilde{x}\tilde{x}}, \\[2mm]
\left(\tilde\rho\left(\tilde e+\frac{\tilde
u^2}{2}\right)\right)_t+\left(\tilde\rho \tilde u\left(\tilde
e+\frac{\tilde u^2}{2}\right)+\tilde p\tilde
u\right)_{\tilde{x}}=\kappa\tilde
\theta_{\tilde{x}\tilde{x}}+(\mu\tilde u\tilde
u_{\tilde{x}})_{\tilde{x}},
\\[4mm]
(\tilde \rho,\tilde u,\tilde \theta)|_{t=0}=(\tilde\rho_0,\tilde
u_0,\tilde\theta_0)(\tilde x)\to(\rho_{\pm},0,\theta_{\pm})
\quad\mbox{as}\quad \tilde{x}\to\pm\infty,
\end{array}\right.
\label{1.1}
\end{equation}
where $\tilde\rho$, $\tilde u$ and $\tilde\theta$ are the density,
the velocity and the absolute temperature, respectively, while
$\mu>0$ is the viscosity coefficient and $\kappa>0$ is the
heat-conductivity coefficients, respectively. It is assumed
throughout the paper that  $\rho_\pm$ and $\theta_\pm$ are
prescribed positive constants and $\R_+=(0,+\infty)$,
$\R_-=(-\infty,0)$ with $\|(\tilde\rho_{0}-\rho_{\pm},\tilde
u_{0},\tilde\theta_{0}-\theta_{\pm})\|_{L^2(\R_{\pm})}$ suitably
small but permitting $|\rho_+-\rho_-|$, $|\theta_+-\theta_-|$ and
$\|(\tilde\rho_{0x},\tilde u_{0x},\tilde\theta_{0x})\|_{L^2(\R)}$
not small. We shall focus our interests on the case of viscous
polytropic ideal gases, so that, the pressure $\tilde p=\tilde
p(\tilde\rho,\tilde\theta)$ and the internal energy $\tilde e=\tilde
e(\tilde\rho,\tilde\theta)$ are related by the second law of
thermodynamics:
\begin{equation}\label{1.2}
\tilde p=R\tilde\rho\tilde\theta,\quad \tilde
e=\frac{R}{\gamma-1}\tilde\theta+const.,
\end{equation}
where $\gamma>1$ is the adiabatic exponent and $R>0$ is the gas
constant.

The purpose is to prove the solvability and stability of the problem
(\ref{1.1}) at $t\in[0,+\infty)$. As there is a local theorem of
existence \cite{SAV} and its references, the main difficulty in
studying the problem in the whole is related to obtaining the a
priori estimate, the constants in which depend only on the
coefficients and the initial data. In this case the local solution
can be extended onto the whole of the length $[0,+\infty)$. When
deducing global estimates, it is convenient to  transform
(\ref{1.1}) to the problem in the Lagrangian coordinate and then
make use of a coordinate transformation to reduce the initial value
problem (\ref{1.1}) into the following form:
\begin{equation}
\left\{
\begin{array}{ll}
v_t-u_x=0,\quad (x,t)\in\R\times\R_+,\\[2mm]
u_t+\left(\frac{R\theta}{v}\right)_x=\mu\left(\frac{u_x}{v}\right)_x,
\\[4mm]
\frac{R}{\gamma-1}\theta_t+R\frac{\theta}{v}u_x=\kappa\left(\frac{\theta_x}{v}\right)_x+\mu\frac{u_x^2}{v},\\[4mm]
(v,u,\theta)|_{t=0}=(v_0,u_0,\theta_0)\to(v_{\pm},0,\theta_{\pm})\quad\mbox{as}\quad
x\to\infty,
\end{array}
\right.\label{1.3}
\end{equation}
where $v_\pm(=\rho^{-1}_{\pm})$ and $\theta_\pm$ are given positive
constants, $\|(v_{0x},u_{0x},\theta_{0x})\|_{L^2(\R)}$ not small and
$v_0,\ \theta_0>0$. Without loss of generality, we set
$1=\theta_+>\theta_->0$.  Here $v=v(x,t),\ u=u(x,t),\
\theta=\theta(x,t)$ are the specific volume, velocity and
temperature as in (\ref{1.1}).

  Up to the present time, some
deep results have been obtained on the asymptotic stability toward
nonlinear waves, viscous shock profiles and viscous rarefaction, for
quite general perturbation of the Navier-Stokes system (\ref{1.3})
and general systems of viscous strictly hyperbolic conservation laws
(see[2--12,14,16--19]). It was observed in \cite{LX,X}, where the
metastability of contact waves was studied for viscous conservation
laws with artificial viscosity dominates the large-time behavior of
solutions. The nonlinear stability of contact discontinuity for the
(full) compressible Navier-Stokes equations was then investigated in
\cite{HMS,HZ}
 for the free  boundary value problem and \cite{HXY,HMX}for the Cauchy problem.
Recently, some problems are call stability of strong viscous
waves(see \cite{NYZ},\cite{Z2013}, \cite{T2015}). These stability
results are shown with some special conditions. Especially,  zero
dissipation result is shown in\cite{MaSX} and $\gamma\to 1$ in
\cite{NYZ} or \cite{HH}. Base on small oscillation, initial
smallness perturbation or zero dissipation (and so on),
Navier-Stokes equations stability results can be obtained with
special help . However to our best knowledge, there is no any
mathematical literature known for the large-time behaviors of
solutions to the general Cauchy problem (\ref{1.3}) due to various
difficulties. To conquer these difficulties, we find that the
important crucial step is  to improve the time estimates of
$(V_x,U_x,\Theta_x)$ in \cite{TZ2012}. This step make the inequality
of Lemma \ref{4.1} be simpler than \cite{HLM2009,TZ2012}, then when
we use similar skills as \cite{T2015}, we can obtain our uniform
esitmates.

The main purpose of this paper is to justify this unknown problem,
i.e., we will show that for a general initial value of the Cauchy
problem (\ref{1.3}), it is possible to be resolved and stable.
Furthermore, the solution approximate the contact discontinuity
$(v_{\pm},0,\theta_{\pm})$. To deduce the desired stability result
by the elementary  energy method, as describe in \cite{SAV},
\cite{HLM2009} and \cite{TZ2012}, it is sufficient to deduce certain
uniform (with respect to the time variable $t$) energy type
estimates on the solution ($v,u,\theta$) and to establish the upper
boundaries of ($v,v^{-1},\theta,\theta^{-1}$), also the
Poincar$\acute{e}$ type inequality in Lemma \ref{lem4.1} without the
smallness of $|\theta_+-\theta_-|$ is important,
 where the arguments employed in \cite{HMS2003,HMS,HMX,HXY,HLM2009,HH,TZ2012} use both smallness $|\theta_+-\theta_-|$ and $N(t)=\sup_{0\leq \tau\leq t}\|(\varphi,\psi,\zeta)\|_{H^1}$ to overcome these difficulties.

Base on the analysis above, the remainder of this paper is organized
as follows. In section 2, we construct a pair of viscous functions
$(V,U,\Theta)(x,t)$ and check that it is nearly close to
$(v_{\pm},0,\theta_{\pm})$. In section 3, we reformulate the problem
and give the precise statement of our main theorem. Finally, we
complete the proof of the main results by the global a priori
estimates established in Section 4.

Throughout this paper, we shall denote $H^l(\omega)$ the usual
$l-th$ order Sobolev space with the norm
$$\|f\|_{l}=\big(\sum_{j=0}^l\|\partial_x^j f\|^2\big)^{1/2},\ \ \|\cdot\|:=\|\cdot\|_{L^2(\omega)}.$$
For simplicity, we also use $C$ or $C_i$ ($i=1,2,3.....$) to denote
the various positive generic constants; $C(\delta_0)$ or
$C_i(\delta_0)$ ($i=1,2,3.....$) to denote one small constant about
$\delta_0^q$ ($q>0$). And $\partial_x^i=\frac{\partial^i}{\partial
x^i}, \ C_v=\frac{R}{\gamma-1}$.
\section{Preliminaries}

In this section, to study the asymptotic behavior of the solution to
the Cauchy problem (\ref{1.3}), we provide some preliminary lemmas
 that are important for the proof of
Theorem \ref{thm3.1}.

First of all, let $\frac{v_-}{\theta_-}=\frac{v_+}{\theta_+}$ and
\begin{equation}
P(V,\Theta)=R\frac{\Theta}{V}=p_+,\;\;\;U(x,t)=\frac{\kappa(\gamma-1)\Theta_x}{\gamma
R\Theta}. \label{2.1}
\end{equation}
 Here $\Theta(x,t)\ ((x,t)\in\R\times\R_+ )$ is the solution of the following  problem
\begin{equation}
\left\{ \begin{array}{ll}
\Theta_t=a(\ln\Theta)_{xx},\quad a=\frac{\kappa p_+(\gamma-1)}{\gamma R^2}>0,\\[2mm]
 \Theta(x,0)=\Theta_{0}(x)\to \theta_{\pm},\\[2mm]
 \Theta_0(x)=\left(\frac{1}{\sqrt{\pi}}(\theta_+^{1/\delta_0}-\theta_-^{1/\delta_0})\int_0^{\ln(x+\sqrt{1+x^2})}\exp\{-y^2\}dy
 +\frac{\theta_-^{1/\delta_0}+\theta_+^{1/\delta_0}}{2}\right)^{\delta_0}.
\end{array}\right.\label{2.2}
\end{equation} In this paper, we ask $\delta_0$ is a suitably small positive constant and $1/\delta_0$ is an integer.

According to the smallness $\delta_0$, we can find that the
properties of $\Theta_{0}(x)$ can be listed as follows.

 \begin{lem}\label{lem2.1}\begin{eqnarray}\label{2.3} &&\|\Theta_{0x}\|_{L^1(\R)}\leq C\ \ ,
|\Theta_{0x}|\leq C\delta_0,\ \ \|\Theta_{0x}\|^2\leq C \delta_0^2,\
\ \|\Theta_0-\theta_{\pm}\|_{L^1(\R_{\pm})}\leq C,
\nonumber\\
&&\|(\ln\Theta_0)_{xx}\|^2\leq C \delta_0^2,\ \
\|(\ln\Theta_0)_{xxx}\|^2\leq C.
\end{eqnarray}
\end{lem}
\pf In fact, if  $K(x)=\ln(x+\sqrt{1+x^2})$, we can get
\begin{eqnarray}\label{2.4}
\int_{\R_{\pm}}\exp\{-K^2(x)\}dx&=&\int_{\R_{\pm}}\exp\{-K^2(x)\}\frac{\sqrt{1+x^2}}{\sqrt{1+x^2}}dx\nonumber\\
&&=\int_{\R_{\pm}}\exp\{-K^2(x)\}\sqrt{1+x^2}d K(x)\nonumber\\
&&\leq\int_{\R_{\pm}}\exp\{-K^2(x)\}\exp\{|K(x)|\}d K(x)\nonumber\\
&&\leq\int_{\R_{\pm}}\exp\{-K^2(x)+|K(x)|\}d K(x)\leq C.
\end{eqnarray}Set
\begin{equation}\label{2.5}H(x)=\Theta_0^{1/\delta_0}=\frac{\theta_+^{1/\delta_0}+\theta_-^{1/\delta_0}}{2}
+\frac{\theta_+^{1/\delta_0}-\theta_-^{1/\delta_0}}{\sqrt{\pi}}\int_0^{K(x)}\exp\{-x^2\}dx,\end{equation}
from $K_x=(1+x^2)^{-1/2}$ and
\begin{equation}\label{2.6}\Theta_{0x}=\delta_0H^{\delta_0-1}(x)H_x(x)=\delta_0H(x)^{\delta_0-1}\frac{\theta_+^{1/\delta_0}-\theta_-^{1/\delta_0}}{\sqrt{\pi}}K_x(x)
\exp\{-K^2(x)\}>0,\ \ x\in\R,\end{equation} we can get
$\| V_{0x}\|_{L^1(\R)}=\int_{\R}\delta_0H^{\delta_0-1}(x)H_x(x)dx<C$.

 When $x>0$, $K(x)>0$, we can know
\begin{equation}\label{2.7}\frac{\theta_+^{1/\delta_0}-\theta_-^{1/\delta_0}}{H(x)}\leq
C\frac{\theta_+^{1/\delta_0}-\theta_-^{1/\delta_0}}{\theta_+^{1/\delta_0}+\theta_-^{1/\delta_0}}\leq
C .\end{equation}  From (\ref{2.5}), \begin{equation}\label{2.8}\theta_-^{1/\delta_0}\leq H(x)\leq
\theta_+^{1/\delta_0},\end{equation} from (\ref{2.6}) and
$K_x(x)=\frac{1}{\sqrt{1+x^2}}$,
 we can get $
|\Theta_{0x}|\leq C\delta_0$.  Also from (\ref{2.4}),(\ref{2.5}),
(\ref{2.6}) and (\ref{2.7}) we can get
$\|\Theta_{0x}\|_{L^2(\R_+)}\leq C\delta_0$. Similar as above
estimates, it is easy to  check $\|\Theta_{0xx}\|_{L^2(\R_+)}^2\leq
C\delta_0^2$ and $\|\Theta_{0xxx}\|_{L^2(\R_+)}^2\leq C$.

When $\delta_0=\frac{1}{2k+1}$, $k\in \mathbb{N}$ is a suitably
large constant, from the  equality
$a^n-b^n=(a-b)\sum_{i=0}^{n-1}a^{n-1-i}b^i$, $\forall a>0,\ b>0,\ n\in\mathbb{N}_+$ and
(\ref{2.4}),(\ref{2.6}), (\ref{2.7}), we can get that
\begin{eqnarray}\label{2.9}&&\int_{\R_+}|\Theta_0-\theta_+|dx=\int_{\R_+}|H^{\delta_0}-\theta_+|dx\nonumber\\
&&=\int_{\R_+}\frac{|H-\theta_+^{1/\delta_0}|}{\sum_{i=0}^{2k}H^{(2k-i)/(2k+1)}\theta_+^{i/(2k+1)}}dx\nonumber\\
&&\leq
C\int_{\R_+}\frac{(\theta_+^{1/\delta_0}-\theta_-^{1/\delta_0})\exp\{-CK^2(x)\}}{\sum_{i=0}^{2k}H^{(2k-i)/(2k+1)}\theta_+^{i/(2k+1)}}dx\nonumber\\
&&\leq
C(\theta_+^{1/\delta_0}-\theta_-^{1/\delta_0})\sup_{x\in\R_+}H^{1/(2k+1)-1}\leq
C(\theta_++\theta_-).\end{eqnarray}

When $x<0$ , $K(-x)=-K(x)>0$.
 We can obtain from (\ref{2.5}) that
\begin{equation}\label{2.10}H(x)=\Theta_0^{1/\delta_0}=\frac{\theta_+^{1/\delta_0}-\theta_-^{1/\delta_0}}{\sqrt{\pi}}\int_{-\infty}^{K(x)}\exp\{-y^2\}dy+\theta_-^{1/\delta_0}.\end{equation}

From
\begin{eqnarray}\label{2.11}&&\int_{-\infty}^{K(x)}\exp\{-x^2\}dx=\left(\int_{-\infty}^{K(y)}\int_{-\infty}^{K(x)}\exp\{-x^2-y^2\}dxdy\right)^{1/2}\nonumber\\
&&\geq\sqrt{\pi}/2\exp\{-r^2/2\}\Big|_{r=-\infty}^{r=\sqrt{2}|K(x)|}\geq\sqrt{\pi}/2\exp\{-K^2(x)\},\end{eqnarray}
 we can obtain
\begin{equation}
\label{2.12}H(x)^{-1}(\theta_+^{1/\delta_0}-\theta_-^{1/\delta_0})\exp\{-K^2(x)\}\leq
C.\end{equation} Combine with (\ref{2.4}), (\ref{2.6}),(\ref{2.8}) and
(\ref{2.12}), we can get
\begin{equation}
\label{2.13}0<\Theta_{0x}\leq
C\delta_0\pi^{-1/2}K_x(x)H^{\delta_0},\end{equation} then
$$|\Theta_{0x}|\leq C\delta_0,\ \ \int_{\R_-}|\Theta_{0x}|^2dx\leq C\delta_0^2.$$
Also $$\|\Theta_{0xx}\|^2_{L^2(\R_-)}\leq C\delta_0^2,\
\|\Theta_{0xxx}\|^2_{L^2(\R_-)}\leq C.$$ Because
$1=\theta_+>\theta_->0$,  similar as (\ref{2.9})we can get
\begin{eqnarray*}&&\int_{\R_-}|\Theta_0-\theta_-|dx=\int_{\R_+}|H^{\delta_0}-\theta_-|dx\\
&&=\int_{\R_-}\frac{|H-\theta_-^{1/\delta_0}|}{\sum_{i=0}^{2k}H^{(2k-i)/(2k+1)}\theta_-^{i/(2k+1)}}dx\\
&&\leq
C\int_{\R_-}\frac{(\theta_+^{1/\delta_0}-\theta_-^{1/\delta_0})\exp\{-3/4K^2(x)\}}{\sum_{i=0}^{2k}H^{(2k-i)/(2k+1)}\theta_-^{i/(2k+1)}}dx\\
&&\leq
C(\theta_+^{1/\delta_0}-\theta_-^{1/\delta_0})\exp\{-1/2K^2(x)\}\sup_{x\in\R_-}(H\theta_-)^{1/2\delta_0-1/2}.\end{eqnarray*}
Combine with (\ref{2.12}) we can get
$\int_{\R_-}|\Theta_0-\theta_-|dx\leq
C\sqrt{\theta_+^{1/\delta_0}-\theta_-^{1/\delta_0}}\leq C.$

So we can get
$$\int_{\R_-}|\Theta_0-\theta_-|dx\leq C.$$
We finish this lemma. $\Box$

In summary, from (\ref{2.1}) and (\ref{2.2}) we have constructed a pair of functions $(V,U,\Theta)$
satisfies
\begin{equation}
\left\{
\begin{array}{ll}
R\frac{\Theta}{V}=p_+,\\[4mm]
V_t=U_x,\\[2mm]
U_t+(R\Theta/V)_x=\mu\left(\frac{U_x}{V}\right)_x+F,\\[4mm]
\frac{R}{\gamma-1}\Theta_t+R\frac{\Theta}{V}U_x=\kappa\left(\frac{\Theta_x}{V}\right)_x+\mu\frac{U_x^2}{V}+G,\\[4mm]
(V,U,\Theta)(x,0)=(V_0,U_0,\Theta_0)=(\frac{R}{p_+}\Theta_0,\frac{k(\gamma-1)}{\gamma
R}\frac{\Theta_{0x}}{\Theta_0},\Theta_0)\to
(v_{\pm},0,\theta_{\pm}),\ as\ \ x\to\infty.
\end{array}\right.\label{2.14}
\end{equation}
where
\begin{eqnarray}&G(x,t)&=-\mu\frac{U_x^2}{V}=O((\ln\Theta)^2_{xx}),\nonumber\\
&F(x,t)&=\frac{\kappa(\gamma-1)}{\gamma
R}\left\{(\ln\Theta)_{xt}-\mu\left(\frac{(\ln\Theta)_{xx}}{V}\right)_x\right\}\nonumber\\
&&=\frac{k a(\gamma-1)-\mu p_+\gamma}{R
\gamma}\left(\frac{(\ln\Theta)_{xx}}{\Theta}\right)_x.\label{2.15}\end{eqnarray}

Furthermore, from (\ref{2.1}) and (\ref{2.14}), we obtain
\begin{eqnarray}\label{2.16}
(|V_x|+U)\leq C|\Theta_x|,\ |\Theta_x|^2\leq
C\|(\ln\Theta)_x\|\|(\ln\Theta)_{xx}\|,\ |U_x|^2\leq
C\|(\ln\Theta)_{xx}\|\|(\ln\Theta)_{xxx}\|.
\end{eqnarray}
When the time $t\to\infty$, it is easily check that $(V,U,\Theta)$
is nearly close to $(v_{\pm},0,\theta_{\pm})$ . We can proof this result by the
following lemma.

\begin{lem}\label{lem2.2}
\begin{eqnarray}
&&\|(\ln\Theta)_x\|^2+\int_0^t\ \|(\ln\Theta)_{xx}\|^2\
dt\leq
C\delta_0^2.\label{2.17}\\
&&_{(see (\ref{2.27})-(\ref{2.28}))}\nonumber\\
 &&\int_0^t\ \|(\ln\Theta)_{x}\|^2\ dt\leq
C(1+t)^{1/3}.\label{2.18}\\
&&_{{(see (\ref{2.29})-(\ref{2.30}))}}\nonumber\\
 && \|(\ln\Theta)_x\|^2\leq
C(1+t)^{-2/3}.\label{2.19}\\&&_{{(see
(\ref{2.31})-(\ref{2.33}))}}\nonumber\\
&&\|(\ln\Theta)_{xx}\|^2\leq C(1+t)^{-5/3}.\label{2.20}\\&&_{(see
(\ref{2.34})-(\ref{2.37}))}\nonumber\\
&&\|(\ln\Theta)_{xx}\|^2(1+t)+\int_0^t\|\partial^3_x\ln\Theta\|^2(1+t)\
dt\leq
C\delta_0^2.\label{2.21}\\&&_{(see(\ref{2.38}))}\nonumber\\
&& \|\partial^3_x\ln\Theta\|^2\leq
C(1+t)^{-8/3}.\label{2.22}\\&&_{(see
(\ref{2.39})-(\ref{2.40}))}\nonumber\\
&&\|\Theta-\theta_{\pm}\|^2_{L^{\infty}(\R_{\pm})}\leq
C\delta_0^{1/4}(1+t)^{-1/24}.\label{2.23}\\&&_{(see
(\ref{2.42})-(\ref{2.44}))}\nonumber
\end{eqnarray}
\end{lem}
\pf Set
\begin{eqnarray*}&&\theta_2(x,t)=\int_{-\infty}^{+\infty}(4\pi
at)^{-1/2}\Theta_0(h)\exp\{-\frac{(h-x)^2}{4at}\}\
dh,\end{eqnarray*}

\begin{eqnarray}\label{2.24}
 &&\theta_{2t}=a\theta_{2xx},\nonumber\\
 &&\theta_2(x,0)=\Theta_0(x)\to\theta_{\pm},
 \end{eqnarray}
and \begin{eqnarray}\label{2.25} \theta_{2x}
&=&\int_{-\infty}^{+\infty}(4\pi
at)^{-1/2}\Theta_0(z)\exp\{-\frac{(z-x)^2}{4at}\}\frac{z-x}{2at}dz\nonumber\\
 &=&\int_{-\infty}^{+\infty}(4\pi
at)^{-1/2}\Theta_{0z}(z)\exp\{\frac{-(z-x)^2}{4at}\}dz\nonumber\\
&=&\int_{-\infty}^{+\infty}(4\pi
at)^{-1/2}\left(\Theta_{0}(z)-\Theta_0(x)\right)\exp\{\frac{-(z-x)^2}{4at}\}\frac{z-x}{2at}dz.
\end{eqnarray}
 By using H$\ddot{o}$lder inequality , Fubini Theorem, (\ref{2.25}) and
$\|\Theta_{0}-\theta_{\pm}\|_{L^1(\R_{\pm})}<C$,
$\|\theta_{2x}\|\leq C(\delta_0)$ , we can get that
\begin{eqnarray}\label{2.26}&&\int_0^t\int_{-\infty}^{\infty}\theta_{2x}^2dxdt\nonumber\\
&&\leq C\int_1^t\int_{-\infty}^{\infty}(4\pi
at)^{-1}\int_{\R_{\pm}}\left|\Theta_{0}(z)-\theta_{\pm}\right|\exp\{-\frac{(z-x)^2}{4at}\}\frac{|z-x|}{2a\sqrt{t}}t^{-1/2}dz
dxdt\nonumber\\&&+C\int_1^t\int_{-\infty}^{\infty}(4\pi
at)^{-1}\int_{\R_{\pm}}\left|\theta_{\pm}-\Theta_0(x)\right|\exp\{-\frac{(z-x)^2}{4at}\}\frac{|z-x|}{2a\sqrt{t}}t^{-1/2}dz dxdt\nonumber\\
&&\quad+\int_0^1\|\theta_{2x}\|^2d\tau\nonumber\\
 &&\leq C\ln(1+t)+C\leq
C(1+t)^{1/3}+C\leq C(1+t)^{1/3} .\end{eqnarray}

 Now, let's consider the estimates about $\partial_x^i\Theta$ $(i=1,2,3)$ of
 (\ref{2.2}).
 In fact from (\ref{2.2}) we can get
\begin{equation*}(\ln\Theta)_t=a\frac{(\ln\Theta)_{xx}}{\Theta}.\end{equation*}
Both side of it multiply by $(\ln\Theta)_{xx}$ and integrate it with
respect to $\R\times (0,t)$ we can get

\begin{eqnarray}&&\|(\ln\Theta)_x\|^2+\int_0^t\ \|(\ln\Theta)_{xx}\|^2\
dt\leq C\|(\ln\Theta_{0})_x\|^2.\label{2.27}\end{eqnarray} When
combine with (\ref{2.3}),
 we can get
\begin{equation}\|(\ln\Theta)_x\|^2+\int_0^t\ \|(\ln\Theta)_{xx}\|^2\ dt\leq
C\delta_0^2.\label{2.28}\end{equation}

On the other hand , integrate
$\left((\ref{2.2})_1-(\ref{2.24})_1)\right)\times(\Theta-\theta_2)$
 in $\R\times (0,t)$ and combine with Cauchy-Schwarz
inequality, we can get
\begin{eqnarray}&&\|\Theta-\theta_2\|^2+\int_0^t\|(\ln\Theta)_x\|^2\
dt \leq C\int_0^t\|{\theta_2}_x\|^2\ dt.\label{2.29}\end{eqnarray}
Insert (\ref{2.26}) and (\ref{2.27}) to (\ref{2.29}) we can get
\begin{equation}\int_0^t\ \|(\ln\Theta)_{x}\|^2\ dt\leq
C(1+t)^{1/3}. \label{2.30}\end{equation} That is (\ref{2.18}).

 Next, integrate $(\ref{2.2})_1\times
\Theta^{-1}(\ln\Theta)_{xx}(1+t)$ in $\R\times(0,t)$, we can get
\begin{eqnarray}0=a\int_0^t\int_{\R}\ \frac{(\ln\Theta)^2_{xx}}{\Theta}(1+t)\
dxdt+\int_0^t\int_{\R}\ \left((\ln\Theta)^2_x\right)_t(1+t)\
dxdt.\label{2.31}\end{eqnarray}  So
\begin{eqnarray}
&&(1+t)\|(\ln\Theta)_x\|^2+\int_0^t\
\int_{\R}(1+t)(\ln\Theta)^2_{xx}\ dx\ dt
dt\nonumber\\
 &&\leq C\|\Theta_{0x}\|^2+\int_0^t\ \int_{\R}\ (\ln\Theta)^2_x\ dx\ dt.\label{2.32}\end{eqnarray}
Combine with (\ref{2.30}) we can get
\begin{eqnarray}
&&(1+t)\|(\ln\Theta)_x\|^2+\int_0^t\
\int_{\R}(1+t)(\ln\Theta)^2_{xx}\ dx\ dt\leq
C(1+t)^{1/3}.\label{2.33}\end{eqnarray} That means
$\|(\ln\Theta)_x\|^2\leq C(1+t)^{-2/3}$, which is(\ref{2.19}).

Again from (\ref{2.2})$_1$ we can get
\begin{equation}(\ln\Theta)_{xt}=a\left(\frac{(\ln\Theta)_{xx}}{\Theta}\right)_x.\label{2.34}
\end{equation}
Both side of (\ref{2.34})multiply by $\partial^3_x\ln\Theta$ and get
\begin{equation}\left((\ln\Theta)_{xt}\partial^2_x(\ln\Theta)\right)_x-1/2(\partial^2_x\ln\Theta)_t
=a\left(\frac{(\ln\Theta)_{xx}}{\Theta}\right)_x\partial^3_x(\ln\Theta).\label{2.35}\end{equation}
 Both side of (\ref{2.35}) multiply by $(1+t)^2$ , integrate it
 with respect to
$\R\times(0,t)$ and combine with Cauchy-Schwarz inequality and
(\ref{2.33})  we have
\begin{eqnarray}&&\|(\ln\Theta)_{xx}\|^2(1+t)^2+\int_0^t\ \int_{\R}\
(1+t)^2(\ln\Theta)_{xxx}^2\ dx\ dt\leq
C(1+t)^{1/3},\label{2.36}\end{eqnarray} which  means
\begin{equation}\|(\ln\Theta)_{xx}\|^2\leq
C(1+t)^{-5/3}.\label{2.37}\end{equation} So we finish (\ref{2.20}).

If both side of (\ref{2.35}) multiply by $(1+t)$, similar as the
proof of (\ref{2.36}),  when combine with (\ref{2.37}) we can get
\begin{equation}\|(\ln\Theta)_{xx}\|^2(1+t)+\int_0^t\ \int_{\R}\ (1+t)(\partial^3_x\ln\Theta)^2\ dx\
dt\leq C\delta_0^2,\label{2.38}\end{equation} which means
(\ref{2.21}).

From (\ref{2.34}) we can get
\begin{equation}
\partial_t(\ln\Theta)_{xx}=a\partial^2_x\left(\frac{(\ln\Theta)_{xx}}{\Theta}\right).\label{2.39}
\end{equation}
When both side of (\ref{2.39}) multiply by
$(\partial^4_x\ln\Theta)(1+t)^3$, then integrate it with respect to
$\R\times(0,t)$, we can get that
 there exists  a constant $\epsilon>0$ such that (\ref{2.39}) can be change to
\begin{eqnarray*}
&&\|\partial^3_x\ln\Theta\|^2(1+t)^3+C\int_0^t\
(1+t)^3\|\partial^4_x\ln\Theta\|^2\ dt\\
&&\leq C+C\int_0^t\ \int_{\R}\
(\partial^3_x\ln\Theta)^2(\ln\Theta)_x^2(1+t)^3\ dx\ dt+C\int_0^t\
\int_{\R}\ (\ln\Theta)_{xx}^4(1+t)^3\ dx\ dt\\
&&\quad+C\int_0^t\ \int_{\R}\
(\partial^2_x\ln\Theta)^2(\ln\Theta)_x^4(1+t)^3\ dx\ dt+C\int_0^t\
\int_{\R}\ (\partial_x^3\ln\Theta)^2(1+t)^2\ dx\ dt\\
&&\leq C\int_0^t\
\|(\ln\Theta)_x\|^2\|\partial^3_x\ln\Theta\|\|\partial^4_x\ln\Theta\|(1+t)^3\
dt+C\int_0^t\
\|(\ln\Theta)_{xx}\|^3\|\partial^3_x\ln\Theta\|(1+t)^3\ dt\\
&&\quad+\int_0^t\ \|(\ln\Theta)_{xx}\|^4\|(\ln\Theta)_x\|^2(1+t)^3\
dt+C(1+t)^{1/3}+C\\
&&\leq \epsilon\int_0^t\ \|\partial^4_x\ln\Theta\|^2(1+t)^3\
dt+C/\epsilon\int_0^t\ \|\partial^3_x\ln\Theta\|^2(1+t)^2\ dt\\
&&\quad+C/\epsilon \int_0^t\ \|\partial^2_x\ln\Theta\|^2(1+t)\
dt+C(1+t)^{1/3}.
\end{eqnarray*}
By using (\ref{2.33}) and (\ref{2.36}) we can get
\begin{equation}\|\partial^3_x\ln\Theta\|^2(1+t)^3+\int_0^t\
(1+t)^3\|\partial^4_x\ln\Theta\|^2\ dt\leq
C(1+t)^{1/3}.\label{2.40}\end{equation} This means (\ref{2.22})
finished.

 When we change $(1+t)^3$  to
$(1+t)^2$ and combine with (\ref{2.21}), we can get

\begin{equation}\|\partial^3_x\ln\Theta\|^2(1+t)^2+\int_0^t\
(1+t)^2\|\partial^4_x\ln\Theta\|^2\ dt\leq
C.\label{2.41}\end{equation}

From (\ref{2.2}) we can
get\begin{equation*}(\Theta-\Theta_0)_t(\Theta-\Theta_0)=a\left((\ln\Theta)_x(\Theta-\Theta_0)\right)_x-a(\ln\Theta)_x(\Theta-\Theta_0)_x.\end{equation*}
When integrate both sides of it in $\R\times[0,t]$, we can get
\begin{eqnarray}\|\Theta-\Theta_0\|^2\leq
C\|\Theta_{0x}\|_{L^1}\int_0^t\|\Theta_x\|_{L^{\infty}}d\tau\leq
C\|\Theta_{0x}\|_{L^1}\int_0^t\|\Theta_x\|^{1/2}\|\Theta_{xx}\|^{1/2}d\tau.\label{2.42}\end{eqnarray}
From  (\ref{2.42}), Lemma \ref{lem2.1}, (\ref{2.19}) and
(\ref{2.20}) we can obtain
\begin{equation}\|\Theta-\Theta_0\|^2\leq C(1+t)^{5/12}.\label{2.43}\end{equation}
From Lemma \ref{lem2.1}, (\ref{2.43}),  (\ref{2.17}) and
(\ref{2.19})
 we can get
\begin{eqnarray}\|\Theta-\theta_{\pm}\|^2_{L^{\infty}(\R_{\pm})}&\leq
&C\|\Theta-\theta_{\pm}\|_{L^2(\R_{\pm})}\|\Theta_x\|^{3/4}\|\Theta_x\|^{1/4}\nonumber\\
&\leq&
C\left(\|\Theta-\Theta_0\|_{L^2(\R_{\pm})}+\|\Theta_0-\theta_{\pm}\|_{L^2(\R_{\pm})}\right)\|\Theta_x\|^{3/4}\|\Theta_x\|^{1/4}\nonumber\\
&\leq& C\delta_0^{1/4}(1+t)^{-1/24}.\label{2.44}\end{eqnarray}

 So we finish this lemma.$\Box$

We can obtain from $|(V-v_{\pm},U,\Theta-\theta_{\pm})|^2(x,t)\leq
C\|(V-v_{\pm},U,\Theta-\theta_{\pm})\|\|(V_x,U_x,\Theta_x)\|$,
(\ref{2.1}),(\ref{2.16}) and Lemma \ref{lem2.2} that for
$x\in\R_{\pm}$,
\begin{equation}\lim_{t\to\infty}|(V,U,\Theta)|(x,t)=(v_{\pm},0,\theta_{\pm}).\label{2.45}\end{equation}
If
$$\|(v-V,u-U,\theta-\Theta)\|_{L^{\infty}(\R)}\to 0,\ \ t\to\infty,$$
we can get that the asymptotic stability results to $(v,u,\theta)$
is $(v_{\pm},0,\theta_{\pm})$. We remain this stability proof at the
end of this paper.

\section{Reformulation and main result}

Now, let $(v,u,\theta)$ be the solution to the problem  (\ref{1.3})
and let $(V,U,\Theta)$ be the solution to (\ref{2.14}). Denote
\begin{eqnarray}
&&\varphi(x,t)=v(x,t)-V(x,t),\nonumber\\
&&\psi(x,t)=u(x,t)-U(x,t),\nonumber\\
&&\zeta(x,t)=\theta(x,t)-\Theta(x,t).\label{3.1}
\end{eqnarray}
Combining (\ref{2.14}) and (\ref{1.3}), the original problem can be
reformulated as
\begin{equation}\left\{
    \begin{array}{lll}
      \varphi_t=\psi_x, &  \\
      \psi_t-(\frac{R\Theta}{vV}\varphi)_x+(\frac{R\zeta}{v})_x=-\mu(\frac{U_x}{vV}\varphi)_x+\mu(\frac{\psi_x}{v})_x-F, & \\
      \frac{R}{\gamma-1}\zeta_t+\frac{R\theta}{v}(\psi_x+U_x)-\frac{R\Theta}{V}U_x
=k(\frac{\zeta_x}{v})_x-k(\frac{\Theta_x\varphi}{vV})_x+\mu(\frac{{u_x}^2}{v}-\frac{{U_x}^2}{V})-G,
&
    \end{array}
  \right.
\label{3.2}\end{equation} and
 \begin{eqnarray}\label{3.3}
(\varphi_0,\psi_0,\zeta_0)=\left(v(x,0)-V(x,0),u(x,0)-U(x,0),\theta(x,0)-\Theta(x,0)\right)
\end{eqnarray}
$$(\varphi_0,\zeta_0)(x)\in H_0^1(0,\infty),\ \ \psi_0(x)\in H^1(0,\infty).$$

From (\ref{2.14}), it is easy to check that the initial  values in
(\ref{3.3}) satisfy
$$(\varphi,\psi,\zeta)(x,0)\to(0,0,0),\ \ \mathrm{as}\ \ x\to \pm\infty.$$ Moreover, for an interval $I\in [0,\infty)$ , we define
the function space
$$X(I)=\left\{(\varphi,\psi,\zeta)\in C(I,H^1)|\varphi_x\in
L^2(I;L^2), (\psi_x,\zeta_x)\in L^2(I;H^1)\right\}.$$ Our  main
results of this paper now reads as follows.
\begin{thm}\label{thm3.1}If
  $(v_0-v_{\pm},u_0,\theta_0-\theta_{\pm})\in H^2(\R_{\pm})\cap L^1(\R_{\pm})$
  ,  $\|(\varphi_0,\psi_0,\zeta_0)\|$ is suitably
  small,
$\frac{v_-}{\theta_-}=\frac{v_+}{\theta_+}$ and
$|\theta_+-\theta_-|$ not small,
  (\ref{3.2}) has a global solution $(\varphi,\psi,\zeta)$  satisfying
$(\varphi,\psi,\zeta)\in X([0,\infty))$, and when $t\to\infty$,
$$\|(\varphi,\psi,\zeta)\|_{L^{\infty}(\R_{\pm})}\to 0.$$
\end{thm}

We prove Theorem \ref{thm3.1} by combining the local existence and
the global-in-time a priori estimates. The local existence of the
solution is well known (e.g., see \cite{SAV,HMS} ), so we omit it
here for brevity. To prove the global existence part of Theorem
\ref{thm3.1}, the same as the asymptotic stability result, we need
to establish the following a priori estimate.
\begin{pro}\label{pro3.1}{\rm(A priori Estimate)} Let $(\varphi,\psi,\zeta)(x,t)\in
X([0,T])$ be a solution of problem (\ref{3.2}) for a constant $T>0$.
Set $C$ is a positive constant only depending on $C_v,\ R,\ \mu,\
\theta_{\pm},\ v_{\pm}$ and $\|(\varphi_0,\psi_0,\zeta_0)\|_1$. If
$\|(\varphi_0,\psi_0,\zeta_0)\|_{L^2(\R)}$ is suitably small and
$$1<\bar{N}_1(T)=\sup_{t\in[0,T]}\{m_v^{-1},M_{v},m_{\theta}^{-1},M_{\theta},\|(\varphi,\psi,\zeta)\|_1\}\leq 2\left(C\|(\varphi_0,\psi_0,\zeta_0)\|_1^2+C+1\right)^{1/2}$$
with $0<m_{v}=v^{-1}(x,t)\leq v(x,t)\leq M_{v}$, $0<m_{\theta}\leq
\theta(x,t)\leq M_{\theta}$, then when $t\in[0,+\infty)$,
$(\varphi,\psi,\zeta)(t)$ satisfies
\begin{eqnarray}\label{3.4}
\sup_{t\in[0,+\infty)}
\|(\varphi,\psi,\zeta)\|_1^2+\int_0^t\left\{\|\varphi_x\|^2+\|(\psi_x,\zeta_x)\|_1^2\right\}d\tau\leq
C\|(\varphi_0,\psi_0,\zeta_0)\|_1^2+C. \end{eqnarray}
\end{pro}

This proposition means  that  for any $T>0$, all the properties of
$X([0,T])$ have uniform bounds, so the  solution's time interval
$[0,T]$ can be extend onto $[0,+\infty)$.

\section{Proof of Theorem \ref{thm3.1} }
Before establishing (\ref{3.4}), we must obtain the upper and lower
boundaries of $v(x,t)$ and $\theta(x,t)$. Here, we set the initial
data
of (\ref{1.3}) are sufficiently smooth functions and set \begin{eqnarray}\label{4.1}&&\Phi(z)=z-\ln z-1,\nonumber\\
&&\Psi(z)=z^{-1}+\ln z-1,\end{eqnarray} where $\Phi'(1)=\Phi(1)=0$
is a strictly convex function around $z=1$. Similar to the proof in
\cite{SAV} or \cite{HMS}, we can get

\begin{lem}\label{lem4.1} If  $C$ is a positive constant independenting of $x$ and $t$, when $\delta_0$ is a small constant,
we can get
\begin{eqnarray*}&&\int_0^t\int_{\R}\Theta_x^2(\varphi^2+\zeta^2)dxd\tau\leq
C(\delta_0)\int_0^t\|(\varphi_x,\zeta_x)\|^2d\tau+C(\delta_0).
\end{eqnarray*}
\end{lem}
\pf \begin{eqnarray*}
&&\int_0^t\int_{\R}\Theta_x^2(\zeta^2+\varphi^2)dxd\tau\\
&&\leq\int_0^t\int_{\R}\Theta_x^2(\|\zeta\|\|\zeta_x\|+\|\varphi\|\|\varphi_x\|)dxd\tau\\
&&\leq C\int_0^t
(\|\zeta_x\|+\|\varphi_x\|)^2\|\Theta_x\|^{1/4}d\tau+C\int_0^t\|\Theta_x\|^{15/4}d\tau.
\end{eqnarray*}
From (\ref{2.17}) and (\ref{2.19}) we can get \begin{eqnarray*}
\int_0^t\int_{\R}\Theta_x^2(\zeta^2+\varphi^2)dxd\tau\leq
C(\delta_0)\int_0^t (\|\zeta_x\|+\|\varphi_x\|)^2d\tau+C(\delta_0).
\end{eqnarray*} That we finish this lemma. $\Box$

\begin{lem}\label{lem4.2}
 If  $C$ is a positive constant independenting of $x$ and $t$, when $\delta_0$ is a small constant,
we can get
\begin{eqnarray*}
&&\int_{\R}\left(R\Theta\Phi\left(\frac{v}{V}\right)+\frac{1}{2}\psi^2+C_v\Theta\Phi\left(\frac{\theta}{\Theta}\right)\right)dx+\int_0^t\
\|(\psi_x/\sqrt{v\theta},\zeta_x/(\theta\sqrt{v}))\|^2\
d\tau\\
&&\leq C(\delta_0)+C(\delta_0)\int_0^t\ \|\varphi_x\|^2\
d\tau+C\|(\varphi_0,\psi_0,\zeta_0)\|^2.
\end{eqnarray*}
\end{lem}
\pf  Similar to the proof in \cite{SAV,HMS} and use the definition
of (\ref{4.1}),
 we deduce from (\ref{3.2}) that
\begin{eqnarray*}
&&\left(\frac{\psi^2}{2}+R\Theta\Phi\left(\frac{v}{V}\right)+C_v\Theta\Phi\left(\frac{\theta}{\Theta}\right)\right)_t\nonumber\\
&&+\mu\frac{\Theta\psi_x^2}{v\theta}+\kappa\frac{\Theta\zeta_x^2}{v\theta^2}+H_x+Q
=\mu\left(\frac{\psi\psi_x}{v}\right)_x-F\psi-\frac{\zeta
G}{\theta},
\end{eqnarray*}
where
$$H=R\frac{\zeta\psi}{v}-R\frac{\Theta\varphi\psi}{vV}+\mu\frac{U_x\varphi\psi}{vV}-\kappa\frac{\zeta\zeta_x}{v\theta}+\kappa\frac{\Theta_x\varphi\zeta}{v\theta V},$$
and
\begin{eqnarray*}
Q&=&p_+\Phi\left(\frac{V}{v}\right)U_x+\frac{p_+}{\gamma-1}\Phi\left(\frac{\Theta}{\theta}\right)U_x-\frac{\zeta}{\theta}(p_+-p)U_x-\mu\frac{U_x\varphi\psi_x}{vV}\\
&&-\kappa
\frac{\Theta_x}{v\theta^2}\zeta\zeta_x-\kappa\frac{\Theta\Theta_x}{v\theta^2V}\varphi\zeta_x-2\mu\frac{U_x}{v\theta}\zeta\psi_x
+\kappa\frac{\Theta_x^2}{v\theta^2V}\varphi\zeta+\mu\frac{U_x^2}{v\theta
V}\varphi\zeta\\
&:=&\sum_{i=1}^9Q_i.
\end{eqnarray*}
Note that $p=R\theta/v$, $p_+=R\Theta/V$, combine h the definition
of $U$, $\bar{N}_1$ with (\ref{2.1}), use integration by parts and
Cauchy-Schwarz inequality, we can get

\begin{eqnarray}\label{4.2}
Q_1+Q_2&=&Ra\left(\Phi\left(\frac{V}{v}\right)(\ln\Theta)_x\right)_x+\frac{Ra}{\gamma-1}\left(\Phi\left(\frac{\Theta}{\theta}\right)(\ln\Theta)_x\right)_x\nonumber\\
&&-aR(\ln\Theta)_x\left(\frac{V\varphi_x\varphi-V_x\varphi^2}{Vv^2}\right)\nonumber\\
&&-a\frac{p_+}{\gamma-1}(\ln\Theta)_x\left(\frac{\Theta\zeta_x\zeta-\Theta_x\zeta^2}{\Theta\theta^2}\right)\nonumber\\
&&\geq\left(p_+\Phi\left(\frac{V}{v}\right)U+\frac{p_+}{\gamma-1}\Phi\left(\frac{\Theta}{\theta}\right)U\right)_x\nonumber\\
&&-C^{1/2}(\delta_0)(\frac{\zeta_x^2}{v\theta^2}+\varphi_x^2)-C^{-1/2}(\delta_0)\bar{N}_1^4\Theta_x^2(\zeta^2+\varphi^2).
\end{eqnarray}
Using $p-p_+=\frac{R\zeta-p_+\varphi}{v}$ and the definition of
$\bar{N}_1$, $U$ and $V$, we can get
\begin{equation*}
Q_3\geq\frac{R\zeta-p_+\varphi}{v}(\frac{\zeta}{\theta}U_x)\geq\left(\frac{R\zeta^2U}{v\theta}-\frac{p_+\zeta\varphi
U}{\theta
v}\right)_x-C^{1/2}(\delta_0)(\frac{\zeta_x^2}{v\theta^2}+\varphi_x^2)-C^{-1/2}(\delta_0)\bar{N}_1^{8}\Theta_x^2(\zeta^2+\varphi^2).
\end{equation*}
And again from the definition of $\bar{N}_1$, Cauchy-Schwarz
inequality we know
\begin{eqnarray*}
(Q_4+Q_7)+(Q_5+Q_6+Q_8)+Q_9&\geq&-C^{-1/2}(\delta_0)\bar{N}_1^4|(\ln\Theta)_{xx}|^2
-C^{1/2}(\delta_0)\frac{\psi_x^2}{v\theta}\nonumber\\
&&-C^{1/2}(\delta_0)\frac{\zeta_x^2}{v\theta^2}-C^{-1/2}(\delta_0)\bar{N}_1^3\Theta_x^2(\zeta^2+\varphi^2).
\end{eqnarray*}
At the end we use the definition of $F$ and $G$ in (\ref{2.15}) then
combine the definition of $\bar{N}_1$ with the general inequality
skills as above to get
\begin{eqnarray}\label{4.3}
-F\psi-G\frac{\zeta}{\theta}&=&-\frac{\kappa a(\gamma-1)-\mu
p_+\gamma}{R\gamma}\left(\frac{(\ln\Theta)_{xx}}{\Theta}\right)_x\psi\nonumber\\
&&+\frac{\mu p_+}{R\Theta}\left(\frac{\kappa
(\gamma-1)}{R\gamma}(\ln\Theta)_{xx}\right)^2\frac{\zeta}{\theta}\nonumber\\
&\leq&-\frac{\kappa a(\gamma-1)-\mu
p_+\gamma}{R\gamma}\left(\frac{(\ln\Theta)_{xx}}{\Theta}\psi\right)_x+\frac{\kappa
a(\gamma-1)-\mu
p_+\gamma}{R\gamma}\frac{(\ln\Theta)_{xx}}{\Theta}\psi_x\nonumber\\
&&+\frac{\mu p_+}{R\Theta}\left(\frac{\kappa
(\gamma-1)}{R\gamma}(\ln\Theta)_{xx}\right)^2\frac{\zeta}{\theta}\nonumber\\
&\leq&-\frac{\kappa a(\gamma-1)-\mu
p_+\gamma}{R\gamma}\left(\frac{(\ln\Theta)_{xx}}{\Theta}\psi\right)_x+C^{1/2}(\delta_0)\frac{\psi_x^2}{v\theta}\nonumber\\
&&+C^{-1/2}(\delta_0)\bar{N}_1^2(\ln\Theta)_{xx}^2.
\end{eqnarray}
Integrating each term from (\ref{4.2}) to (\ref{4.3})) in
$\R\times(0,t)$ , using Lemma \ref{lem4.1} and Lemma \ref{lem2.2} at
the end,  we get that for a small $C(\delta_0)>0$,
\begin{eqnarray*}
&&\int_{\R}\left(R\Theta\Phi\left(\frac{v}{V}\right)+\frac{1}{2}\psi^2+C_v\Theta\Phi\left(\frac{\theta}{\Theta}\right)\right)dx+\int_0^t\
\|(\psi_x/\sqrt{v\theta},\zeta_x/\sqrt{v\theta^2})\|^2\
d\tau\\
&&\leq C(\delta_0)(1+\bar{N}_1^{8})+C(\delta_0)\int_0^t\
\|\varphi_x\|^2\ d\tau+C\|(\varphi_0,\psi_0,\zeta_0)\|^2.
\end{eqnarray*}
Then we finish this lemma. $\Box$

\begin{lem}\label{lem4.3} If $\|(\varphi_0,\psi_0,\zeta_0)\|$ and $C(\delta_0)>0$ are suitably small and $\|\varphi_{0x}\|$ not small,  we can get
\begin{eqnarray*}
\|\varphi_x/v\|^2+\int_0^t\|\varphi_x/v\sqrt{R\theta/v}\|^2
d\tau\leq
C\|(\varphi_0,\varphi_{0x},\psi_0,\zeta_0)\|^2+C(\delta_0),
\end{eqnarray*}with $C$ independent of $x$ and $t$.
\end{lem}
\pf Set $\bar{v}=\frac{v}{V}$, take it into (\ref{3.2})$_1$,
(\ref{3.2})$_2$  ($p=R \theta/v$) to get
$$
\psi_t+p_x=\mu\left (\frac{\bar{v}_x}{\bar{v}}\right)_t-F,
$$
Both sides of last equation multiply by $\bar v_x/\bar v$ to get
\begin{eqnarray}
&&\left(\frac{\mu}{2}\left (\frac{\bar v_x}{\bar
v}\right)^2-\psi\frac{\bar v_x}{\bar v}\right)_t
+\frac{R\theta}{v}\left (\frac{\bar v_x}{\bar
v}\right)^2+\left(\psi\frac{\bar v_t}{\bar v}\right)_x\nonumber\\
&&\quad=\frac{\psi_x^2}{v}+U_x\left
(\frac{1}{v}-\frac{1}{V}\right)\psi_x+\frac{R\zeta_x}{v}\frac{\bar
v_x}{\bar v}-\frac{R\theta}{v}\left
(\frac{1}{\Theta}-\frac{1}{\theta}\right)\Theta_x\frac{\bar
v_x}{\bar v}+F\frac{\bar v_x}{\bar v}.\label{4.4}
\end{eqnarray}

On the other hand if we integrate (\ref{4.4}) in $R\times(0,t)$, the
right side of it is less than
\begin{eqnarray*}&&C\left(\int_0^t\ \|(\zeta_x/\sqrt{v\theta},\psi_x/\sqrt{v})\|^2\
d\tau+\bar{N}_1^2\int_0^t\ \int_{\R}\ \Theta_x^2\zeta^2\ dx\
d\tau\right)\\&&+C\bar{N}_1\int_0^t\int_{\R}U_x^2\varphi^2
dxd\tau+C\bar{N}_1^2\int_0^t\int_{\R}|F|^2dxd\tau+1/4\int_0^t\int_{\R}\frac{R\theta}{v}\left(\frac{\bar{v}_x}{\bar{v}}\right)^2\
dxd\tau.\end{eqnarray*} Use Lemma \ref{lem4.1} to the term
$\int_0^t\ \int_{\R}\ \Theta_x^2(\varphi^2+\zeta^2)\ dxd\tau$ and
combine with the definition of $(V,U,\Theta)$,
$\int_0^t\int_{\R}(\ref{4.4})dxd\tau$ can be change to
\begin{eqnarray*}&&\int_0^t\|\sqrt{\frac{R\theta}{v}}\frac{\bar{v}_x}{\bar{v}}\|^2\
d\tau+\|\frac{\bar{v}_x}{\bar{v}}\|^2-C\|\psi\|^2-C\|\psi_0\|^2-C\int_{\R}\frac{\bar{v}_x}{\bar{v}}(x,0)^2\
dx\nonumber\\ &&\leq C\left(\int_0^t\
\|(\zeta_x/\sqrt{v\theta},\psi_x/\sqrt{v})\|^2\
d\tau+C(\delta_0)\bar{N}_1^2\int_0^t\ \int_{\R}\ \zeta_x^2\ dx\
d\tau+C(\delta_0)\bar{N}_1^2\right)+C\|\varphi_{0x}\|^2.
\end{eqnarray*}
Insert $C_1(\frac{\varphi_x}{v})^2-C\bar{N}_1^2\varphi^2V_x^2\leq
(\frac{\bar{v}_x}{\bar{v}})^2\leq
C_3(\frac{\varphi_x}{v})^2+C\bar{N}_1^2\varphi^2V_x^2$  to last
inequality and combine with Lemma \ref{lem4.2}, Lemma \ref{lem2.2}
and the definition of $\bar{N}_1$ we can obtain that
\begin{eqnarray*}&&\int_0^t\|\varphi_x/v\sqrt{R\theta/v}\|^2\
d\tau+\|\varphi_x/v\|^2\nonumber\\&&\leq
C\|(\varphi_0,\varphi_{0x},\psi_0,\zeta_0)\|^2+C(\delta_0)(1+\bar{N}_1^8)
+C\int_0^t\ \|(\zeta_x/\sqrt{v\theta},\psi_x/\sqrt{v})\|^2\ d\tau.
\end{eqnarray*}
 So we finish this lemma. $\Box$

 When $C(\delta_0)$ and $\|(\varphi_0,\psi_0,\zeta_0)\|$ are suitably small and use it to control the terms
about $\bar{N}_1$, insert Lemma \ref{lem4.3} to Lemma
\ref{lem4.2}, we obtain the
first a priori energy estimate
\begin{eqnarray}\label{4.5}\sup_{t\in[0,+\infty)}\int_{\R}(\varphi^2+\psi^2+\zeta^2)dx+\int_0^t\
\|(\psi_x,\zeta_x)\|^2\ d\tau&\leq
C\|(\varphi_0,\psi_0,\zeta_0)\|^2+C(\delta_0).\end{eqnarray}

\begin{lem}\label{lem4.4} For a suitably small constant
$C(\delta_0)>0$ and $\|(\psi_{0x},\zeta_{0x})\|$ not small, the
solution of (\ref{3.2}) satisfies
\begin{eqnarray*}
\|(\psi_x,\zeta_x)\|^2+\int_0^t\|(\psi_{xx},\zeta_{xx})(\tau)\|^2d\tau\leq
C\|(\varphi_0,\psi_0,\zeta_0)\|_1^2+C(\delta_0),
\end{eqnarray*}with $C$ independent of $x$ and $t$.
\end{lem}
\pf Multiplying $(\ref{3.2})_2$ and $(\ref{3.2})_3$ by $\psi_{xx}$
and $\zeta_{xx}$, respectively, and summing up the resulting
equations, we find
\begin{eqnarray*}
&&\left(\frac{1}{2}\psi_x^2+{\frac{R}{2(\gamma-1)}}\zeta_x^2\right)_t+\mu\frac{\psi_{xx}^2}{v}+\kappa\frac{\zeta_{xx}^2}{v}\\
&&\quad= (\varphi_x+V_x)\left(\mu\frac{\psi_x
\psi_{xx}}{v^2}+\kappa\frac{\zeta_x\zeta_{xx}}{v^2}\right)+\left(\mu\frac{U_x\varphi}{v
V}-R\frac{\Theta\varphi}{v
V}+R\frac{\zeta}{v}\right)_x\psi_{xx}\\
&&\qquad+\left\{\kappa\left(\frac{\Theta_x\varphi}{vV}\right)_x
-\mu\left(\frac{u_x^2}{v}-\frac{U_x^2}{V}\right)+R\frac{\theta\psi_x}{v}+R\frac{U_x\zeta}{v}-R\frac{U_x\Theta\varphi}{vV}\right\}\zeta_{xx}
\\[2mm]
&&\qquad+\left(F\psi_{xx}+G\zeta_{xx}\right)+\left(\psi_t\psi_x+\frac{R}{\gamma-1}\zeta_t\zeta_x\right)_x=I_1+I_2+I_3+I_4+I_5.
\end{eqnarray*}
which, integrated over $\R\times (0,t)$, gives
\begin{eqnarray}
&&\|(\psi_x,\zeta_x)(t)\|^2+\int_0^t\|(\psi_{xx},\zeta_{xx})(\tau)\|^2d\tau\nonumber\\
&&\quad\leq
C\|(\psi_{0x},\zeta_{0x})\|^2+C\sum_{i=1}^5\left|\int_0^t\int_{-\infty}^\infty
I_idxd\tau\right|.\label{4.6}
\end{eqnarray}
We now estimate  $\iint I_idxd\tau$ ($i=1,\ldots,5$). First, using
 Cauchy-Schwarz inequality, we
have
\begin{eqnarray*}
&&\int_0^t \int_{-\infty}^{\infty}|I_1| dx d\tau\leq
\epsilon\int_0^t\|(\psi_{xx},\zeta_{xx})\|^2d\tau+
\frac{C}{\epsilon}\int_0^t\int_{-\infty}^\infty(V_x^2+\varphi_x^2)(\psi_x^2+\zeta_x^2)dx
d\tau\nonumber\\
&&\quad\leq\epsilon\int_0^t\|(\psi_{xx},\zeta_{xx})\|^2d\tau+
\frac{C}{\epsilon}\int_0^t\left(\|V_x\|_{L^\infty}^2\|(\psi_x,\zeta_x)\|^2
+\|(\psi_x,\zeta_x)\|_{L^\infty}^2\|\varphi_x\|^2\right)d\tau,
\end{eqnarray*}
which, combined with
$$
\|(\psi_x,\zeta_x)\|_{L^\infty}^2\leq
C\|(\psi_x,\zeta_x)\|\|(\psi_{xx},\zeta_{xx})\|,
$$
yields
\begin{equation}
\int_0^t \int_{-\infty}^{\infty}|I_1| dx
d\tau\leq2\epsilon\int_0^t\|(\psi_{xx},\zeta_{xx})(\tau)\|^2d\tau+
\frac{C}{\epsilon^2}\left(C(\delta_0)+N^4(T)\right)\int_0^t\|(\psi_x,\zeta_x)(\tau)\|^2
d\tau.\label{4.7}
\end{equation}

Noting that
\begin{eqnarray*}
|I_2|&\leq&
C\left(|U_{xx}||\varphi|+|U_x||\varphi_x|+|U_x||\varphi||\varphi_x|+|U_x||V_x||\varphi|+|\Theta_x||\varphi|+|\varphi_x|\right.\\
&&\left.+|\varphi||\varphi_x|+|V_x||\varphi|+|\zeta_x|+|\zeta||\varphi_x|+|V_x||\zeta|\right)|\psi_{xx}|\\
&\leq&C\left(|\varphi_x|+|\zeta_x|+|U_x||\varphi_x|+|U_x||\varphi||\varphi_x|+|\varphi||\varphi_x|+|\zeta||\varphi_x|\right)|\psi_{xx}|\\
&&+C\left(|U_{xx}|+|U_x||V_x|+|\Theta_x|+|V_x|\right)(|\varphi|+|\zeta|)|\psi_{xx}|,
\end{eqnarray*}
using (\ref{4.5}), Lemma \ref{lem4.1},  Cauchy-Schwarz inequalities
and Sobolev inequalities, we infer that
\begin{eqnarray}
\int_0^t\int_{-\infty}^\infty
|I_2|dxd\tau&\leq&\epsilon\int_0^t\|\psi_{xx}(\tau)\|^2d\tau+\frac{C}{\epsilon}\int_0^t\|(\varphi_x,\zeta_x)(\tau)\|^2
d\tau\nonumber\\
&&+\frac{C}{\epsilon}\int_0^t \int_{-\infty}^\infty
\Theta_x^2(\varphi^2+\zeta^2) dx d\tau\nonumber\\
&&+\frac{C}{\epsilon}\int_0^t\
\|(\varphi,\zeta)\|^2\left(\|\partial_x^3(\ln\Theta)\|^2+\|\varphi_x\|^2\right)\
d\tau\nonumber\\
&\leq& \epsilon\int_0^t\ \|\psi_{xx}\|^2\
d\tau+\frac{C}{\epsilon}\left(\|(\varphi_0,\psi_0,\zeta_0)\|^2+\|\varphi_{0x}\|^2+C(\delta_0)\right)
.\label{4.8}
\end{eqnarray}
In a similar manner, we also have
\begin{eqnarray}
\int_0^t\int_{-\infty}^\infty
|I_3|dxd\tau&\leq&\epsilon\int_0^t\|(\psi_{xx},\zeta_{xx})\|(\tau)^2d\tau+\frac{C}{\epsilon}\int_0^t\|(\varphi_x,\psi_x)(\tau)\|^2
d\tau\nonumber\\
&&+\frac{C}{\epsilon}\int_0^t \int_{-\infty}^\infty
\Theta_x^2(\varphi^2+\zeta^2) dx d\tau+\frac{C}{\epsilon}\int_0^t\ \|U_x\|^2\ d\tau\nonumber\\
&\leq&\epsilon\int_0^t\|(\psi_{xx},\zeta_{xx})\|(\tau)^2d\tau+\frac{C}{\epsilon}
\left(\|(\varphi_0,\psi_0,\theta_0)\|^2+C(\delta_0)\right),\label{4.9}
\end{eqnarray}
where the following inequality is used
\begin{eqnarray*}
\|\psi_x^2\zeta_{xx}\|_{L^1}&\leq&\|\psi_x\|_{L^\infty}\|\psi_x\|\|\zeta_{xx}\|\leq
C\|\psi_x\|^{3/2}_{L^2}\|\psi_{xx}\|^{1/2}\|\zeta_{xx}\|\\
&\leq&
\epsilon\|(\psi_{xx},\zeta_{xx})\|^2+\frac{C}{\epsilon}N^4(T)\|\psi_x\|^2.
\end{eqnarray*}
By virtue of Cauchy-Schwarz inequality
\begin{equation}
\big|\int_0^t\int_{-\infty}^\infty\ I_4\ dxd\tau\big|\leq
\frac{C}{\epsilon}\int_0^t\ \left(\|F\|^2+\|G\|^2\right)\
d\tau+{\epsilon}\int_0^t\
\left(\|\psi_{xx}\|^2+\|\zeta_{xx}\|^2\right)\ d\tau.\label{4.10}
\end{equation}
Noting that
\begin{equation}
\left|\int_0^t\int_{-\infty}^\infty I_5dxd\tau\right|=0,\label{4.11}
\end{equation}
inserting the inequalities from (\ref{4.7}) to (\ref{4.11}) into
(\ref{4.6}) and choosing $\epsilon>0$ suitably small, recalling the
definition of $G,\ F$, combining Lemma \ref{lem4.3}, (\ref{4.5}) and
Lemma \ref{lem2.2}, we conclude that for a
$M_v(t)=\|v(x,t)\|_{L^{\infty}(\R)}$
\begin{eqnarray}
\|(\psi_x,\zeta_x)(t)\|^2+\int_0^t\|(\psi_{xx},\zeta_{xx})(\tau)\|^2d\tau\leq
CM_v(t)\left(\|(\varphi_0,\psi_0,\zeta_0)\|_1^2+1\right)\label{4.12}.
\end{eqnarray}
Let (\ref{4.12}) combine with (\ref{4.5})  and Lemma \ref{lem4.3},
we can obtain that for small $\|(\varphi_0,\psi_0,\zeta_0)\|$ and
$C(\delta_0)$ we can get
$$|v-V|\leq C\|\varphi\|\|\varphi_x\|<C\bar{N}_1^5\left(\|(\varphi_0,\psi_0,\zeta_0)\|+C(\delta_0)\right)$$
 which means there exist positive constant $C_5$ and $C_6$
 independent of $x$ and $t$
such that $C_5<v<C_6$. Now we inserting this bounds of $v$ into
(\ref{4.12}) we can get for a positive constant $C$ independent of
$x$ and $t$ such that
\begin{eqnarray*}
\|(\psi_x,\zeta_x)(t)\|^2+\int_0^t\|(\psi_{xx},\zeta_{xx})(\tau)\|^2d\tau\leq
C\left(\|(\varphi_0,\psi_0,\zeta_0)\|_1^2+1\right).
\end{eqnarray*}
We finish this lemma. $\Box$

 In
all, from (\ref{4.5}) , Lemma \ref{lem4.3} and Lemma \ref{lem4.4} we
finish Proposition \ref{pro3.1}.

To finish Theorem \ref{thm3.1} now let's consider the stability
result.

In fact from Proposition \ref{pro3.1} we can get
$$
\int_0^\infty\left(\left|\frac{d}{dt}\|\psi_x(t)\|^2\right|+\left|\frac{d}{dt}\|\zeta_x(t)\|^2\right|+\left|\frac{d}{dt}\|\varphi_x(t)\|^2\right|\
\right)d\tau\leq C\|(\varphi_0,\psi_0,\zeta_0)\|_1^2+C.
$$

It means
$$
\|(\varphi,\psi,\zeta)(t)\|_{L^\infty}^2\leq
2\|(\varphi,\psi,\zeta)(t)\|\|(\varphi_x,\psi_x,\zeta_x)(t)\|\to0\quad\mbox{when}\quad
t\to\infty.
$$
Combine with (\ref{2.45}) we finish the staility result of the
theorem.

\end{document}